\documentclass[12pt]{article}%
\usepackage{sw20bams}
\usepackage{amsmath}
\usepackage{amsfonts}
\usepackage{amssymb}
\usepackage{graphicx}%
\providecommand{\U}[1]{\protect\rule{.1in}{.1in}}
\begin{document}

\title{Triangles Formed via Poisson Nearest Neighbors}
\author{Steven Finch}
\date{December 21, 2017}
\maketitle

\begin{abstract}
We start with certain joint densities (for sides and for angles) corresponding
to pinned Poissonian triangles in the plane, then discuss analogous results
for staked and anchored triangles. \ 

\end{abstract}

\footnotetext{Copyright \copyright \ 2017 by Steven R. Finch. All rights
reserved.}A planar triangle is \textbf{pinned Poissonian} if one vertex $A$ is
fixed at the origin and the other vertices $B$, $C$ are the two nearest
neighboring particles to $(0,0)$ of a unit intensity Poisson process. \ Hence
$\left\Vert B\right\Vert ^{2}$ and $\left\Vert C\right\Vert ^{2}-\left\Vert
B\right\Vert ^{2}$ are independent Exponential($\pi$) variables, where
$\left\Vert \cdot\right\Vert $ denotes Euclidean norm. \ Let $a$, $b$, $c$
denote the sides opposite the random vertices; clearly $c<b$ and $a<2b$ due to
the triangle inequality. A\ Jacobian determinant calculation given in Appendix
1 yields
\[
\left\{
\begin{array}
[c]{l}%
8\pi\dfrac{x\,y\,z}{\sqrt{(x+y+z)(-x+y+z)(x-y+z)(x+y-z)}}\exp\left(
-\pi\,y^{2}\right) \\
\;\;\;\;\;\;\;\;\;\;\;\;\;\;\;\text{if }y-z<x<y+z,\;\;y>z,\;\;z>0\\
0\;\;\;\;\;\;\;\;\;\;\;\;\;\text{otherwise}%
\end{array}
\right.
\]
as the trivariate density $f(x,y,z)$ for $a$, $b$, $c$. As a consequence, the
univariate density for $c$ is
\[%
\begin{array}
[c]{ccc}%
2\pi x\exp\left(  -\pi\,x^{2}\right)  , &  & x>0,
\end{array}
\]
the univariate density for $b$ is
\[%
\begin{array}
[c]{ccc}%
2\pi^{2}x^{3}\exp\left(  -\pi\,x^{2}\right)  , &  & x>0
\end{array}
\]
and the univariate density for $a$ is
\[%
\begin{array}
[c]{ccc}%
\pi\,x\operatorname*{erfc}\left(  \sqrt{\pi}\,x/2\right)  , &  & x>0.
\end{array}
\]
Densities for $c$ and $b$ were known earlier \cite{SS-Pois, Ha-Pois, SN-Pois},
but that for $a$ seems to be new. Table 1 contains moments for these
variables. \ For example, the cross-correlation coefficient
\[
\rho(a,b)=\frac{\operatorname*{Cov}(a,b)}{\sqrt{\operatorname*{Var}%
(a)\operatorname*{Var}(b)}}=\frac{8}{3}\sqrt{\frac{32-9\pi}{(-64+27\pi)\pi}%
}\approx0.636
\]
indicates strong positive dependency. An exact expression for
$\operatorname*{E}(a\,c)=0.49181215...$ remains open.

\smallskip Table 1 \textit{Pinned Moments of Angles, Sides and some Products
and Ratios}%

\[%
\begin{tabular}
[c]{|l|l|l|}\hline
Variable & Mean & Mean Square\\\hline
$\alpha$ & $\pi/2$ & $\pi^{2}/3$\\\hline
$\beta$ & $\pi/4+1/\pi$ & $1+\pi^{2}/12$\\\hline
$\gamma$ & $\pi/4-1/\pi$ & $-1/2+\pi^{2}/12$\\\hline
$\alpha\beta$ & $1/4+\pi^{2}/12$ & -\\\hline
$\beta\gamma$ & $-1/4+\pi^{2}/12$ & -\\\hline
$\gamma\alpha$ & $-1/4+\pi^{2}/12$ & -\\\hline
$a$ & $8/(3\pi)$ & $3/\pi$\\\hline
$b$ & $3/4$ & $2/\pi$\\\hline
$c$ & $1/2$ & $1/\pi$\\\hline
$ab$ & $64/(9\pi^{2})$ & -\\\hline
$bc$ & $4/(3\pi)$ & -\\\hline
$ca$ & $0.49181215...$ & -\\\hline
$a/b$ & $32/(9\pi)$ & $3/2$\\\hline
$b/a$ & $4/\pi$ & $\infty$\\\hline
$b/c$ & $2$ & $\infty$\\\hline
$c/b$ & $2/3$ & $1/2$\\\hline
$c/a$ & $(1+2G)/\pi$ & $\infty$\\\hline
$a/c$ & $(5+2G)/\pi$ & $\infty$\\\hline
area & $4/(3\pi^{2})$ & $3/(8\pi^{2})$\\\hline
\end{tabular}
\
\]

\smallskip From this, another calculation in Appendix 1 gives
\[
\left\{
\begin{array}
[c]{lll}%
\dfrac{2}{\pi}\dfrac{\sin(x)\sin(x+y)}{\sin(y)^{3}} &  & \text{if }%
0<x<\pi\text{ and }\dfrac{\pi-x}{2}<y<\pi-x,\\
0 &  & \text{otherwise}%
\end{array}
\right.
\]
as the bivariate density for angles $\alpha$, $\beta$ opposite sides $a$, $b$.
Let $\gamma=\pi-\alpha-\beta$ and let $\wp$ denote the probability that a
pinned Poissonian triangle is obtuse. While $\alpha$ is Uniform[$0,\pi$], the
univariate density for $\beta$ is%
\[
\left\{
\begin{array}
[c]{lll}%
\dfrac{1}{2\pi}+\dfrac{1-3\cos(x)^{2}}{2\pi\sin(x)^{2}}+\dfrac{x\cos(x)}%
{\pi\sin(x)^{3}} &  & \text{if }0<x<\pi/2,\\
\dfrac{1}{\pi\sin(x)^{2}}+\dfrac{(\pi-x)\cos(x)}{\pi\sin(x)^{3}} &  & \text{if
}\pi/2<x<\pi
\end{array}
\right.
\]
and the univariate density for $\gamma$ is%
\[
\left\{
\begin{array}
[c]{lll}%
(4/\pi)\cos(x)^{2} &  & \text{if }0<x<\pi/2,\\
0 &  & \text{if }\pi/2<x<\pi.
\end{array}
\right.
\]
\ These imply that
\begin{align*}
\wp &  =\operatorname*{P}(\alpha>\pi/2)+\operatorname*{P}(\beta>\pi
/2)+\operatorname*{P}(\gamma>\pi/2)\\
&  =1/2+1/4+0=3/4
\end{align*}
because a triangle can have at most one obtuse angle. \ See \cite{ES1-Pois,
Fi2-Pois} for alternative approaches for computing $\wp$.

The ratio of a pair of sides is of interest \cite{DS-Pois}. \ We find that the
univariate density for $a/b$ is
\[%
\begin{array}
[c]{ccc}%
\dfrac{2\,x}{\pi}\arccos\left(  \dfrac{x}{2}\right)  , &  & 0<x<2,
\end{array}
\]
the univariate density for $b/a$ is
\[%
\begin{array}
[c]{ccc}%
\dfrac{1}{x^{3}}-\dfrac{2}{\pi\,x^{3}}\arcsin\left(  \dfrac{1}{2\,x}\right)
, &  & x>\dfrac{1}{2},
\end{array}
\]
the univariate density for $b/c$ is
\[%
\begin{array}
[c]{ccc}%
\dfrac{2}{x^{3}}, &  & x>1,
\end{array}
\]
the univariate density for $c/b$ is
\[%
\begin{array}
[c]{ccc}%
2\,x, &  & 0<x<1,
\end{array}
\]
the univariate density for $c/a$ is%
\[
\left\{
\begin{array}
[c]{lll}%
2x\dfrac{1+x^{2}}{\left(  1-x^{2}\right)  ^{3}} &  & \text{if }0<x<1/2,\\
-\dfrac{2\left(  -1+x^{2}\right)  \sqrt{-1+4x^{2}}-\pi x^{2}\left(
1+x^{2}\right)  +6x^{2}\left(  1+x^{2}\right)  \arcsin\left(  1/(2x)\right)
}{\pi x\left(  -1+x^{2}\right)  ^{3}} &  & \text{if }x>1/2
\end{array}
\right.
\]
and the univariate density for $a/c$ is%
\[
\left\{
\begin{array}
[c]{lll}%
-x\dfrac{2x\left(  1-x^{2}\right)  \sqrt{4-x^{2}}-\pi\left(  1+x^{2}\right)
+6\left(  1+x^{2}\right)  \arcsin\left(  x/2\right)  }{\pi\left(
1-x^{2}\right)  ^{3}} &  & \text{if }0<x<2,\\
2x\dfrac{1+x^{2}}{\left(  -1+x^{2}\right)  ^{3}} &  & \text{if }x>2.
\end{array}
\right.
\]
It is remarkable that Catalan's constant \cite{Fi3-Pois}
\[
G=%
{\displaystyle\sum\limits_{k=0}^{\infty}}
\frac{(-1)^{k}}{(2k+1)^{2}}%
\]
appears in expressions for both $\operatorname*{E}(c/a)$ and
$\operatorname*{E}(a/c)$, as well as in other geometric probability settings
\cite{Fi4-Pois, Fi5-Pois}.

A planar triangle is \textbf{staked Poissonian} if one vertex $A$ is fixed at
$(0,0)$, another vertex $B$ is fixed at $(1,0)$ and the third vertex $C$ is
the nearest neighboring particle to $(0,0)$ of a unit intensity Poisson
process. The term \textit{stake} (as in \textquotedblleft staking a
tent\textquotedblright) was only recently introduced in this context
\cite{Fi6-Pois}. \ Let $C=(u,v)$. \ Clearly%
\[%
\begin{array}
[c]{ccc}%
\tan(\alpha)=\dfrac{v}{1-u}, &  & \tan(\beta)=\dfrac{v}{u}.
\end{array}
\]
The Jacobian determinant of the transformation $(u,v)\mapsto(\alpha,\beta)$ is%
\[
|J|=\left\vert
\begin{array}
[c]{cc}%
\dfrac{v}{(1-u)^{2}+v^{2}} & \dfrac{1-u}{(1-u)^{2}+v^{2}}\\
\dfrac{-v}{u^{2}+v^{2}} & \dfrac{u}{u^{2}+v^{2}}%
\end{array}
\right\vert =\dfrac{v}{\left(  u^{2}+v^{2}\right)  \left[  (1-u)^{2}%
+v^{2}\right]  }.
\]
Solving for $u$, $v$ in terms of $\alpha$, $\beta$, we obtain%
\[%
\begin{array}
[c]{ccc}%
u=\dfrac{\tan(\alpha)}{\tan(\alpha)+\tan(\beta)}, &  & v=\dfrac{\tan
(\alpha)\tan(\beta)}{\tan(\alpha)+\tan(\beta)}.
\end{array}
\]
Substituting these expressions into a nonstandard bivariate normal density
\[
\pi\exp\left[  -\pi\left\{  u^{2}+v^{2}\right\}  \right]
\]
and dividing by $|J|$ yields%
\[
\exp\left[  -\pi\frac{\sin(\alpha)^{2}}{\sin(\alpha+\beta)^{2}}\right]
\frac{\sin(\alpha)\sin(\beta)}{\sin(\alpha+\beta)^{3}}.
\]
Multiplying by $2$ gives the correct normalization. \ While $\alpha$ is
Uniform[$0,\pi$], the univariate density for $\beta$ is decidedly not so.
\ The acuteness probability can be found exactly: \ \
\begin{align*}
1-\wp &  =2%
{\displaystyle\int\limits_{0}^{\frac{\pi}{2}}}
{\displaystyle\int\limits_{\frac{\pi}{2}-\alpha}^{\frac{\pi}{2}}}
\exp\left[  -\pi\frac{\sin(\alpha)^{2}}{\sin(\alpha+\beta)^{2}}\right]
\frac{\sin(\alpha)\sin(\beta)}{\sin(\alpha+\beta)^{3}}d\beta\,d\alpha\\
&  =\frac{1}{2}\left[  e^{-\pi/2}I_{0}\left(  \frac{\pi}{2}\right)
-\operatorname*{erfc}\left(  \sqrt{\pi}\right)  \right]  =0.1725524698...
\end{align*}
where $I_{0}$ is the $0^{\text{th}}$ modified Bessel function of the first
kind. \ Densities for sides $a$, $b$ are possible but omitted. \ Such
triangles seem to be mentioned in \cite{DS-Pois, DL-Pois} without further elaboration.

\smallskip Table 2 \textit{Staked Moments of Angles and of a Product}%

\[%
\begin{tabular}
[c]{|l|l|l|}\hline
Variable & Mean & Mean Square\\\hline
$\alpha$ & $\pi/2$ & $\pi^{2}/3$\\\hline
$\beta$ & $0.34306160...$ & $0.20825399...$\\\hline
$\alpha\beta$ & $0.43825535...$ & -\\\hline
\end{tabular}
\ \ \
\]

\smallskip A planar triangle is \textbf{anchored Poissonian} if one vertex $A$
is fixed at $(-1/2,0)$, another vertex $B$ is fixed at $(1/2,0)$ and the third
vertex $C$ is the nearest neighboring particle to $(0,0)$ of a unit intensity
Poisson process. The term \textit{anchoring} (as in \textquotedblleft
anchoring a ship\textquotedblright) again was only recently introduced
\cite{Fi6-Pois}. \ Let $C=(u,v)$. \ Clearly%
\[%
\begin{array}
[c]{ccc}%
\tan(\alpha)=\dfrac{v}{\frac{1}{2}-u}=\dfrac{2v}{1-2u}, &  & \tan
(\beta)=\dfrac{v}{\frac{1}{2}+u}=\dfrac{2v}{1+2u}%
\end{array}
\]
and the corresponding Jacobian determinant is%
\[
|J|=\left\vert
\begin{array}
[c]{cc}%
\dfrac{4v}{(1-2u)^{2}+4v^{2}} & \dfrac{2(1-2u)}{(1-2u)^{2}+4v^{2}}\\
-\dfrac{4v}{(1+2u)^{2}+4v^{2}} & \dfrac{2(1+2u)}{(1+2u)^{2}+4v^{2}}%
\end{array}
\right\vert =\dfrac{16\,v}{\left[  (1-2u)^{2}+4v^{2}\right]  \left[
(1+2u)^{2}+4v^{2}\right]  }.
\]
Solving for $u$, $v$ in terms of $\alpha$, $\beta$, we obtain%
\[%
\begin{array}
[c]{ccc}%
u=\dfrac{1}{2}\dfrac{\tan(\alpha)-\tan(\beta)}{\tan(\alpha)+\tan(\beta)}, &  &
v=\dfrac{\tan(\alpha)\tan(\beta)}{\tan(\alpha)+\tan(\beta)}.
\end{array}
\]
Substituting these expressions into the nonstandard bivariate normal density
and dividing by $|J|$ yields%
\[
\exp\left[  -\frac{\pi}{4}\frac{\sin(\alpha-\beta)^{2}+4\sin(\alpha)^{2}%
\sin(\beta)^{2}}{\sin(\alpha+\beta)^{2}}\right]  \frac{\sin(\alpha)\sin
(\beta)}{\sin(\alpha+\beta)^{3}}.
\]
Multiplying by $2$ gives the correct normalization. \ The univariate densities
for $\alpha$ and $\beta$ are identical but decidedly not uniform. The
acuteness probability can be found exactly: \ \
\begin{align*}
1-\wp &  =2%
{\displaystyle\int\limits_{0}^{\frac{\pi}{2}}}
{\displaystyle\int\limits_{\frac{\pi}{2}-\alpha}^{\frac{\pi}{2}}}
\exp\left[  -\frac{\pi}{4}\frac{\sin(\alpha-\beta)^{2}+4\sin(\alpha)^{2}%
\sin(\beta)^{2}}{\sin(\alpha+\beta)^{2}}\right]  \frac{\sin(\alpha)\sin
(\beta)}{\sin(\alpha+\beta)^{3}}d\beta\,d\alpha\\
&  =e^{-\pi/4}-\operatorname*{erfc}\left(  \sqrt{\pi}/2\right)
=0.2458467223...
\end{align*}
This value is slightly larger for anchored triangles than for staked
triangles. \ Densities for sides $a$, $b$ are possible but again omitted.

\smallskip Table 3 \textit{Anchored Moments of Angles and of a Product}%

\[%
\begin{tabular}
[c]{|l|l|l|}\hline
Variable & Mean & Mean Square\\\hline
$\alpha$ & $0.71706372...$ & $0.92490176...$\\\hline
$\beta$ & $0.71706372...$ & $0.92490176...$\\\hline
$\alpha\beta$ & $0.39837926...$ & -\\\hline
\end{tabular}
\ \
\]

\smallskip Appendix 2 discusses an evidently unrelated class of random
triangles $T$. \ If $\varphi$, $\psi$ are independent Uniform[$0,\pi$]
variables, then%
\[
(\alpha,\beta)=\left\{
\begin{array}
[c]{lll}%
(\varphi,\psi) &  & \text{if }\varphi+\psi<\pi,\\
(\pi-\psi,\pi-\varphi) &  & \text{if }\varphi+\psi>\pi
\end{array}
\right.
\]
are angles of $T$ at vertices $A=(0,0)$, $B=(1,0)$. \ Note that if
$\varphi+\psi>\pi$, then%
\[
(\pi-\psi)+(\pi-\varphi)=\pi+(\pi-\varphi-\psi)<\pi\text{;}%
\]
it follows that $\alpha+\beta<\pi$ always. Upon constructing a line $L_{A}$
emanating from $A$ with slope $\tan(\alpha)$ and a line $L_{B}$ emanating from
$B$ with slope $-\tan(\beta)$, the remaining vertex%
\[
C=\left(  \frac{\tan(\beta)}{\tan(\alpha)+\tan(\beta)},\frac{\tan(\alpha
)\tan(\beta)}{\tan(\alpha)+\tan(\beta)}\right)
\]
is the point $L_{A}\cap L_{B}$ of intersection. \ Such \textquotedblleft
uniform triangles\textquotedblright\ have appeared before in the literature
\cite{Grf-Pois, Mor-Pois, ES2-Pois}, although perhaps not with the same
specificity as \cite{DS-Pois}.

\section{Acknowledgement}

I am grateful to Daryl Daley for suggesting the study of pinned Poissonian
triangles in $\mathbb{R}^{2}$. \ Analogous work for Poisson processes in
$\mathbb{R}^{d}$ for $d>2$ awaits an interested reader!\ Thanks are also due
to Adrian Baddeley and Rolf Turner for writing an R\ package \textit{spatstat}
\cite{BRT-Pois}, which enables testing of numerical predictions in this essay
via simulation \cite{Fi7-Pois}.

\section{Appendix 1}

Revisiting the beginning, let $R_{1}=\left\Vert B\right\Vert ^{2}$ and
$R_{2}=\left\Vert C\right\Vert ^{2}-\left\Vert B\right\Vert ^{2}$. \ Define
$\theta_{1}$ to be the angle between vector $B$ and the horizontal axis;
define $\theta_{2}$ likewise for vector $C$. \ The joint density for
$(R_{1},R_{2},\theta_{1},\theta_{2})$ is%
\[
\left(  \frac{1}{2\pi}\right)  ^{2}\left(  \pi\,e^{-\pi R_{1}}\right)  \left(
\pi\,e^{-\pi R_{2}}\right)  =\frac{1}{4}e^{-\pi(R_{1}+R_{2})}%
\]
where $R_{i}>0$, $0<$ $\theta_{i}<2\pi$ for $i=1,2$. \ We rewrite this density
in terms of sides $b$, $c$. From $R_{1}=c^{2}$, $R_{2}=b^{2}-c^{2}$ emerges a
Jacobian matrix%
\[
\left(
\begin{array}
[c]{cc}%
0 & 2c\\
2b & -2c
\end{array}
\right)
\]
with absolute determinant $4bc$. \ Thus the joint density for $(b,c,\theta
_{1},\theta_{2})$ is%
\[
b\,c\,e^{-\pi\,b^{2}}%
\]
where $0<c<b$. \ As in \cite{Fi1-Pois}, we integrate out $\theta_{1}$ by
letting $\omega=\theta_{1}-\theta_{2}$ and $\alpha=\left\vert \omega
\right\vert $, then adding contributions at $\alpha$ and $2\pi-\alpha$.
\ Omitting details, the joint density for $(\alpha,b,c)$ comes out as%
\[
4\pi b\,c\,e^{-\pi\,b^{2}}%
\]
where $0<\alpha\,<\pi$. \ We now bring $a$ into the density, removing $\alpha
$. Differentiating the Law of Cosines
\[
a^{2}=b^{2}-2\,b\,c\cos(\alpha)+c^{2}%
\]
with respect to $\alpha$, it is clear that
\begin{align*}
2\,a\,da  &  =2\,b\,c\sin(\alpha)\,d\alpha\\
\  &  =\sqrt{(a+b+c)(-a+b+c)(a-b+c)(a+b-c)}\,d\alpha
\end{align*}
by a formula for area, and hence the density becomes%
\begin{align*}
&  4\pi b\,c\,e^{-\pi\,b^{2}}d\alpha\,db\,dc\\
&  =4\pi b\,c\,e^{-\pi\,b^{2}}\frac{2a}{\sqrt{(a+b+c)(-a+b+c)(a-b+c)(a+b-c)}%
}da\,db\,dc
\end{align*}
as was to be shown.

Let $\Delta=(a+b+c)(-a+b+c)(a-b+c)(a+b-c)$. \ The natural transformation
$(\alpha,\beta,c)\mapsto$ $(a,b,c)$ appearing in \cite{Fi1-Pois} has Jacobian
determinant $a\,b$. \ Using the identities
\[%
\begin{array}
[c]{ccccc}%
\dfrac{a}{c}=\dfrac{\sin(\alpha)}{\sin(\alpha+\beta)}, &  & \dfrac{b}%
{c}=\dfrac{\sin(\beta)}{\sin(\alpha+\beta)}, &  & \dfrac{\sqrt{\Delta}}%
{2c^{2}}=\dfrac{\sin(\alpha)\sin(\beta)}{\sin(\alpha+\beta)}%
\end{array}
\]
we have%
\[
\sin(\alpha+\beta)=\frac{c}{a}\frac{c}{b}\dfrac{\sqrt{\Delta}}{2c^{2}}%
=\dfrac{\sqrt{\Delta}}{2a\,b}%
\]
thus the pinned angle density can be rewritten is
\begin{align*}
&  \ 8\pi\frac{a^{2}b^{2}c}{\sqrt{\Delta}}\exp\left[  -\pi\,b^{2}\right] \\
\  &  =8\pi c^{5}\frac{\sin(\alpha)^{2}\sin(\beta)^{2}}{\sin(\alpha+\beta
)^{4}\sqrt{\Delta}}\exp\left[  -\pi c^{2}\frac{\sin(\beta)^{2}}{\sin
(\alpha+\beta)^{2}}\right] \\
\  &  =4\pi c^{3}\frac{\sin(\alpha)\sin(\beta)}{\sin(\alpha+\beta)^{3}}%
\exp\left[  -\pi c^{2}\frac{\sin(\beta)^{2}}{\sin(\alpha+\beta)^{2}}\right]  .
\end{align*}
Integrating out $c$ is facilitated by observing that
\[%
{\displaystyle\int\limits_{0}^{\infty}}
c^{3}\exp\left(  -\pi\,c^{2}r\right)  dc=\frac{1}{2\pi^{2}r^{2}}%
\]
for $r>0$, therefore the density for $(\alpha,\beta)$ is
\[
\frac{2}{\pi}\frac{\sin(\alpha)\sin(\beta)}{\sin(\alpha+\beta)^{3}}\left(
\frac{\sin(\alpha+\beta)^{2}}{\sin(\beta)^{2}}\right)  ^{2}=\frac{2}{\pi}%
\frac{\sin(\alpha)\sin(\alpha+\beta)}{\sin(\beta)^{3}}.
\]
The restriction $\beta>(\pi-\alpha)/2$ is implied by $b>c$, equivalently,
$\sin(\beta)>\sin(\alpha+\beta)$.

We note that the joint density for $(a,b)$ is%
\[
4\pi\,a\,b\exp\left(  -\pi\,b^{2}\right)  \arccos\left(  \frac{a}{2b}\right)
\]
for $0<a<2b$ and the joint density for $(b,c)$ is%
\[
4\pi^{2}b\,c\exp\left(  -\pi\,b^{2}\right)
\]
for $0<c<b$. \ No closed-form expression of the joint density $f$ for $(a,c)$
is known. \ If $0<a<2c$, then $a-c<c$ and%
\[
f(a,c)=8\pi\,a\,c%
{\displaystyle\int\limits_{c}^{a+c}}
\frac{b\exp\left(  -\pi\,b^{2}\right)  }{\sqrt{\left[  (a+c)^{2}-b^{2}\right]
\left[  b^{2}-(a-c)^{2}\right]  }}db;
\]
if $0<2c<a$, then $c<a-c$ and%
\[
f(a,c)=8\pi\,a\,c%
{\displaystyle\int\limits_{a-c}^{a+c}}
\frac{b\exp\left(  -\pi\,b^{2}\right)  }{\sqrt{\left[  (a+c)^{2}-b^{2}\right]
\left[  b^{2}-(a-c)^{2}\right]  }}db.
\]
The density $g$ for the ratio $z=a/c$ is \cite{Pa-Pois}%
\begin{align*}
g(z)  &  =%
{\displaystyle\int\limits_{0}^{\infty}}
c\,f(z\,c,c)dc\\
&  =\left\{
\begin{array}
[c]{lll}%
8\pi\,z%
{\displaystyle\int\limits_{0}^{\infty}}
c^{3}%
{\displaystyle\int\limits_{c}^{(z+1)c}}
\dfrac{b\exp\left(  -\pi\,b^{2}\right)  }{\sqrt{\left[  (z+1)^{2}c^{2}%
-b^{2}\right]  \left[  b^{2}-(z-1)^{2}c^{2}\right]  }}db\,dc &  & \text{if
}0<z<2,\\
8\pi\,z%
{\displaystyle\int\limits_{0}^{\infty}}
c^{3}%
{\displaystyle\int\limits_{(z-1)c}^{(z+1)c}}
\dfrac{b\exp\left(  -\pi\,b^{2}\right)  }{\sqrt{\left[  (z+1)^{2}c^{2}%
-b^{2}\right]  \left[  b^{2}-(z-1)^{2}c^{2}\right]  }}db\,dc &  & \text{if
}z>2
\end{array}
\right. \\
&  =\left\{
\begin{array}
[c]{lll}%
8\pi\,z%
{\displaystyle\int\limits_{0}^{\infty}}
b\exp\left(  -\pi\,b^{2}\right)
{\displaystyle\int\limits_{\frac{b}{z+1}}^{b}}
\dfrac{c^{3}}{\sqrt{\left[  (z+1)^{2}c^{2}-b^{2}\right]  \left[
b^{2}-(z-1)^{2}c^{2}\right]  }}dc\,db &  & \text{if }0<z<2,\\
8\pi\,z%
{\displaystyle\int\limits_{0}^{\infty}}
b\exp\left(  -\pi\,b^{2}\right)
{\displaystyle\int\limits_{\frac{b}{z+1}}^{\frac{b}{z-1}}}
\dfrac{c^{3}}{\sqrt{\left[  (z+1)^{2}c^{2}-b^{2}\right]  \left[
b^{2}-(z-1)^{2}c^{2}\right]  }}dc\,db &  & \text{if }z>2
\end{array}
\right.
\end{align*}
and the density $h$ for the ratio $w=c/a$ is%
\begin{align*}
h(w)  &  =%
{\displaystyle\int\limits_{0}^{\infty}}
a\,f(a,w\,a)da\\
&  =\left\{
\begin{array}
[c]{lll}%
8\pi\,w%
{\displaystyle\int\limits_{0}^{\infty}}
a^{3}%
{\displaystyle\int\limits_{w\,a}^{(w+1)a}}
\dfrac{b\exp\left(  -\pi\,b^{2}\right)  }{\sqrt{\left[  (w+1)^{2}a^{2}%
-b^{2}\right]  \left[  b^{2}-(w-1)^{2}a^{2}\right]  }}db\,da &  & \text{if
}w>1/2,\\
8\pi\,w%
{\displaystyle\int\limits_{0}^{\infty}}
a^{3}%
{\displaystyle\int\limits_{(1-w)a}^{(w+1)a}}
\dfrac{b\exp\left(  -\pi\,b^{2}\right)  }{\sqrt{\left[  (w+1)^{2}a^{2}%
-b^{2}\right]  \left[  b^{2}-(w-1)^{2}a^{2}\right]  }}db\,da &  & \text{if
}0<w<1/2
\end{array}
\right. \\
&  =\left\{
\begin{array}
[c]{lll}%
8\pi\,w%
{\displaystyle\int\limits_{0}^{\infty}}
b\exp\left(  -\pi\,b^{2}\right)
{\displaystyle\int\limits_{\frac{b}{w+1}}^{\frac{b}{w}}}
\dfrac{a^{3}}{\sqrt{\left[  (w+1)^{2}a^{2}-b^{2}\right]  \left[
b^{2}-(w-1)^{2}a^{2}\right]  }}da\,db &  & \text{if }w>1/2,\\
8\pi\,w%
{\displaystyle\int\limits_{0}^{\infty}}
b\exp\left(  -\pi\,b^{2}\right)
{\displaystyle\int\limits_{\frac{b}{w+1}}^{\frac{b}{1-w}}}
\dfrac{a^{3}}{\sqrt{\left[  (w+1)^{2}a^{2}-b^{2}\right]  \left[
b^{2}-(w-1)^{2}a^{2}\right]  }}da\,db &  & \text{if if }0<w<1/2.
\end{array}
\right.
\end{align*}
Both $g(z)$ and $h(w)$ are readily evaluated. \ The circumstances are less
advantageous in the following:%
\[
\operatorname*{E}(a\,c)=%
{\displaystyle\int\limits_{0}^{\infty}}
{\displaystyle\int\limits_{0}^{2c}}
a\,c\,f(a,c)da\,dc+%
{\displaystyle\int\limits_{0}^{\infty}}
{\displaystyle\int\limits_{0}^{\frac{a}{2}}}
a\,c\,f(a,c)dc\,da
\]
for which the first integral becomes%
\begin{align*}
&
{\displaystyle\int\limits_{0}^{\infty}}
{\displaystyle\int\limits_{0}^{2c}}
8\pi\,a^{2}c^{2}%
{\displaystyle\int\limits_{c}^{a+c}}
\frac{b\exp\left(  -\pi\,b^{2}\right)  }{\sqrt{\left[  (a+c)^{2}-b^{2}\right]
\left[  b^{2}-(a-c)^{2}\right]  }}db\,da\,dc\\
&  =8\pi%
{\displaystyle\int\limits_{0}^{\infty}}
c^{2}%
{\displaystyle\int\limits_{c}^{3c}}
b\exp\left(  -\pi\,b^{2}\right)
{\displaystyle\int\limits_{b-c}^{2c}}
\frac{a^{2}}{\sqrt{\left[  (a+c)^{2}-b^{2}\right]  \left[  b^{2}%
-(a-c)^{2}\right]  }}da\,db\,dc\\
&  =8\pi%
{\displaystyle\int\limits_{0}^{\infty}}
b\exp\left(  -\pi\,b^{2}\right)
{\displaystyle\int\limits_{\frac{b}{3}}^{b}}
c^{2}%
{\displaystyle\int\limits_{b-c}^{2c}}
\frac{a^{2}}{\sqrt{\left[  (a+c)^{2}-b^{2}\right]  \left[  b^{2}%
-(a-c)^{2}\right]  }}da\,dc\,db
\end{align*}
and the second integral becomes%
\begin{align*}
&
{\displaystyle\int\limits_{0}^{\infty}}
{\displaystyle\int\limits_{0}^{\frac{a}{2}}}
8\pi\,a^{2}c^{2}%
{\displaystyle\int\limits_{a-c}^{a+c}}
\frac{b\exp\left(  -\pi\,b^{2}\right)  }{\sqrt{\left[  (a+c)^{2}-b^{2}\right]
\left[  b^{2}-(a-c)^{2}\right]  }}db\,dc\,da\\
&  =8\pi%
{\displaystyle\int\limits_{0}^{\infty}}
a^{2}%
{\displaystyle\int\limits_{\frac{a}{2}}^{\frac{3a}{2}}}
b\exp\left(  -\pi\,b^{2}\right)
{\displaystyle\int\limits_{\left\vert b-a\right\vert }^{\frac{a}{2}}}
\frac{c^{2}}{\sqrt{\left[  (a+c)^{2}-b^{2}\right]  \left[  b^{2}%
-(a-c)^{2}\right]  }}dc\,db\,da\\
&  =8\pi%
{\displaystyle\int\limits_{0}^{\infty}}
b\exp\left(  -\pi\,b^{2}\right)
{\displaystyle\int\limits_{\frac{2b}{3}}^{2b}}
a^{2}%
{\displaystyle\int\limits_{\left\vert b-a\right\vert }^{\frac{a}{2}}}
\frac{c^{2}}{\sqrt{\left[  (a+c)^{2}-b^{2}\right]  \left[  b^{2}%
-(a-c)^{2}\right]  }}dc\,da\,db.
\end{align*}
Both triple integrals can be reduced to double integrals, each involving an
incomplete elliptic integral of the third kind, but further simplification
does not seem to be feasible.

\section{Appendix 2}

Starting from the joint density for angles in $T$:
\[
\left\{
\begin{array}
[c]{lll}%
2/\pi^{2} &  & \text{if }0<\alpha<\pi\text{, }0<\beta<\pi\text{ and }%
\alpha+\beta<\pi,\\
0 &  & \text{otherwise}%
\end{array}
\right.
\]
we find the joint density for sides%
\[
k(a,b)=\left\{
\begin{array}
[c]{lll}%
2/(\pi^{2}a\,b) &  & \text{if }\left\vert 1-a\right\vert <b<1+a\text{ and
}a>0,\\
0 &  & \text{otherwise}%
\end{array}
\right.
\]
which is true because $c=1$ and since the natural transformation
$(\alpha,\beta)\mapsto$ $(a,b)$ has Jacobian determinant $a\,b$. \ Integrating
out $b$, the marginal density for $a$ is%
\[%
\begin{array}
[c]{ccc}%
k_{a}(a)=\dfrac{2}{\pi^{2}}\dfrac{\ln(1+a)-\ln\left\vert 1-a\right\vert }%
{a}, &  & a>0\text{.}%
\end{array}
\]
Note the singularity at $a=1$. \ The density for $z=a/b$ is%
\[%
{\displaystyle\int\limits_{0}^{\infty}}
b\,k(z\,b,b)db=%
{\displaystyle\int\limits_{\frac{1}{1+z}}^{\frac{1}{\left\vert 1-z\right\vert
}}}
\frac{2}{\pi^{2}z\,b}db=\dfrac{2}{\pi^{2}}\dfrac{\ln(1+z)-\ln\left\vert
1-z\right\vert }{z}%
\]
for $z>0$, which interestingly is the same as that for $a$, $b$ and $b/a$ as
well! \ Similar general formulas \cite{Pa-Pois} are applicable to
$x=\max\{a,b\}$:%
\[%
{\displaystyle\int\limits_{0}^{x}}
k(x,b)db+%
{\displaystyle\int\limits_{0}^{x}}
k(a,x)da=\dfrac{4}{\pi^{2}}\dfrac{\ln(x)-\ln\left\vert 1-x\right\vert }{x}%
\]
for $x>1/2$ and $y=\min\{a,b\}$:
\[
k_{a}(y)+k_{b}(y)-%
{\displaystyle\int\limits_{0}^{y}}
k(y,b)db+%
{\displaystyle\int\limits_{0}^{y}}
k(a,y)da=\left\{
\begin{array}
[c]{lll}%
\dfrac{4}{\pi^{2}}\dfrac{\ln(1+y)-\ln(1-y)}{y} &  & \text{if }0<y<1/2,\\
\dfrac{4}{\pi^{2}}\dfrac{\ln(1+y)-\ln(y)}{y} &  & \text{if }y>1/2.
\end{array}
\right.
\]
Our proofs are simpler than those in \cite{DS-Pois}. \ We have not examined,
however, the complicated density for the area of $T$.

\section{Appendix 3}

Here, for completeness' sake, are R\ simulation output results (histograms in
blue) graphed against density expressions found in this paper (curves in red). \

%

\begin{figure}[ptb]%
\centering
\includegraphics[
height=7.9477in,
width=5.1747in
]%
{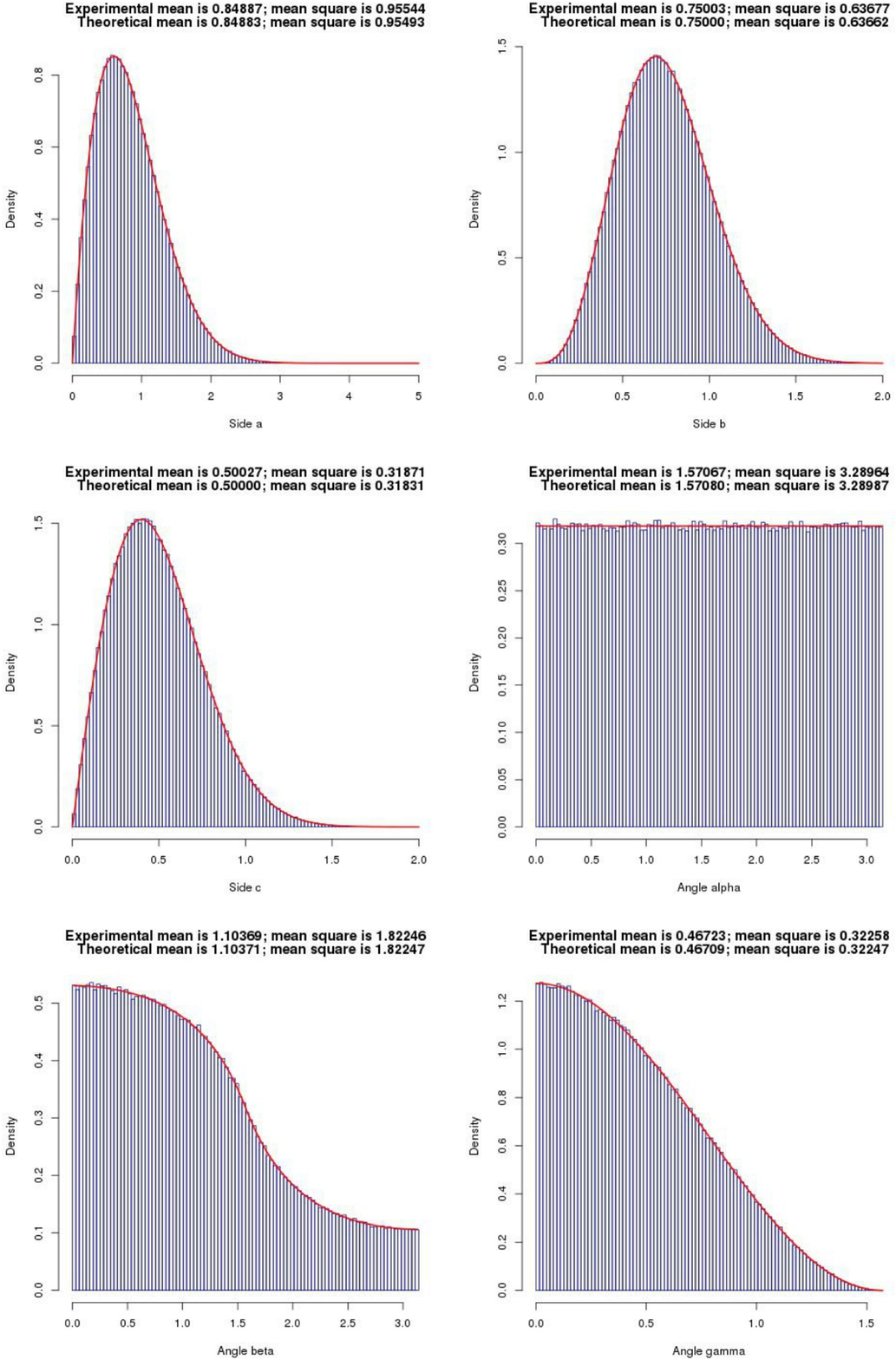}%
\caption{Pinned densities for Poissonian triangle sides and angles.}%
\end{figure}

%

\begin{figure}[ptb]%
\centering
\includegraphics[
height=8.3478in,
width=5.5301in
]%
{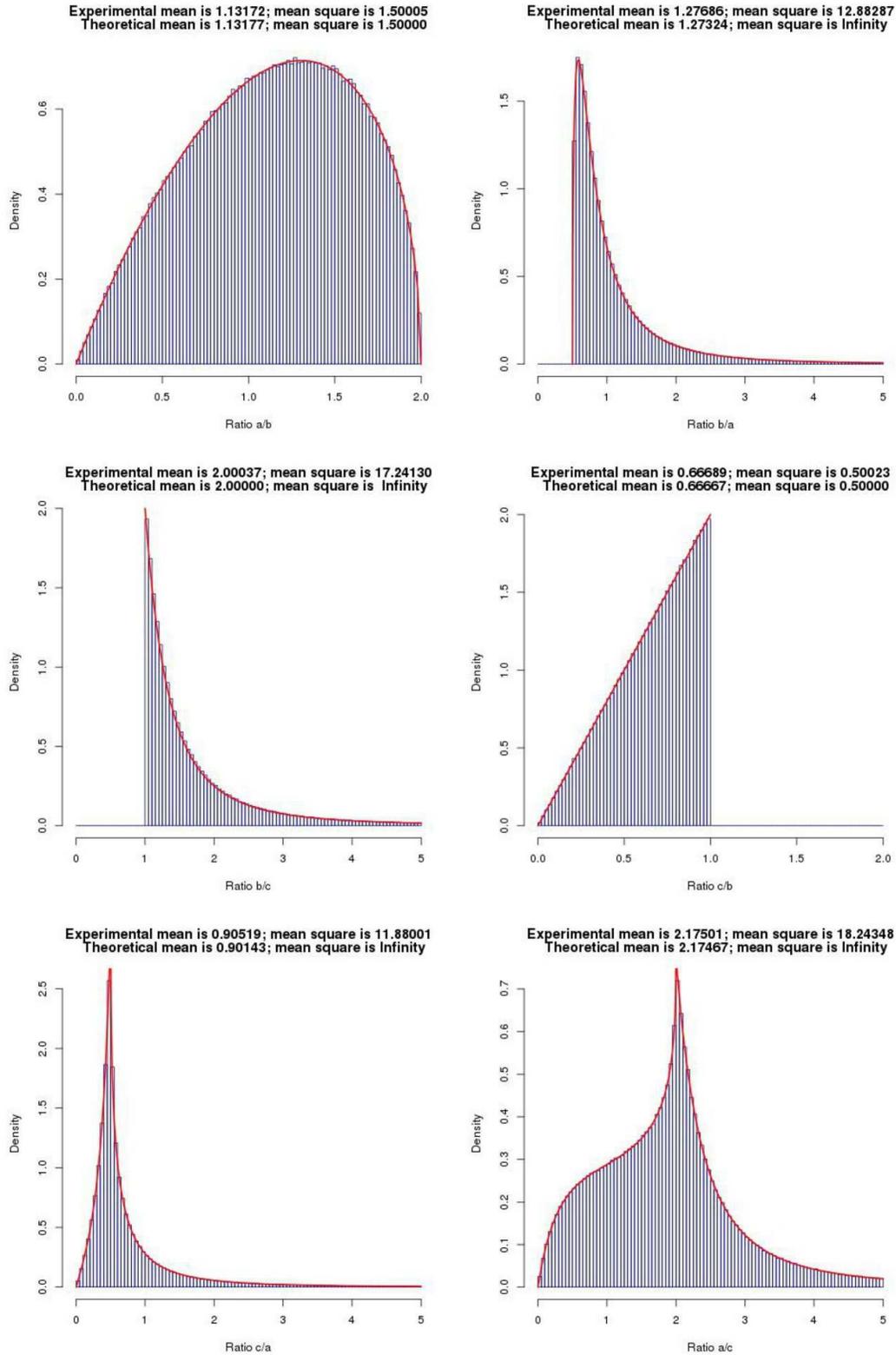}%
\caption{Pinned densities for Poissonian triangle side ratios.}%
\end{figure}

%

\begin{figure}[ptb]%
\centering
\includegraphics[
height=6.3362in,
width=6.3362in
]%
{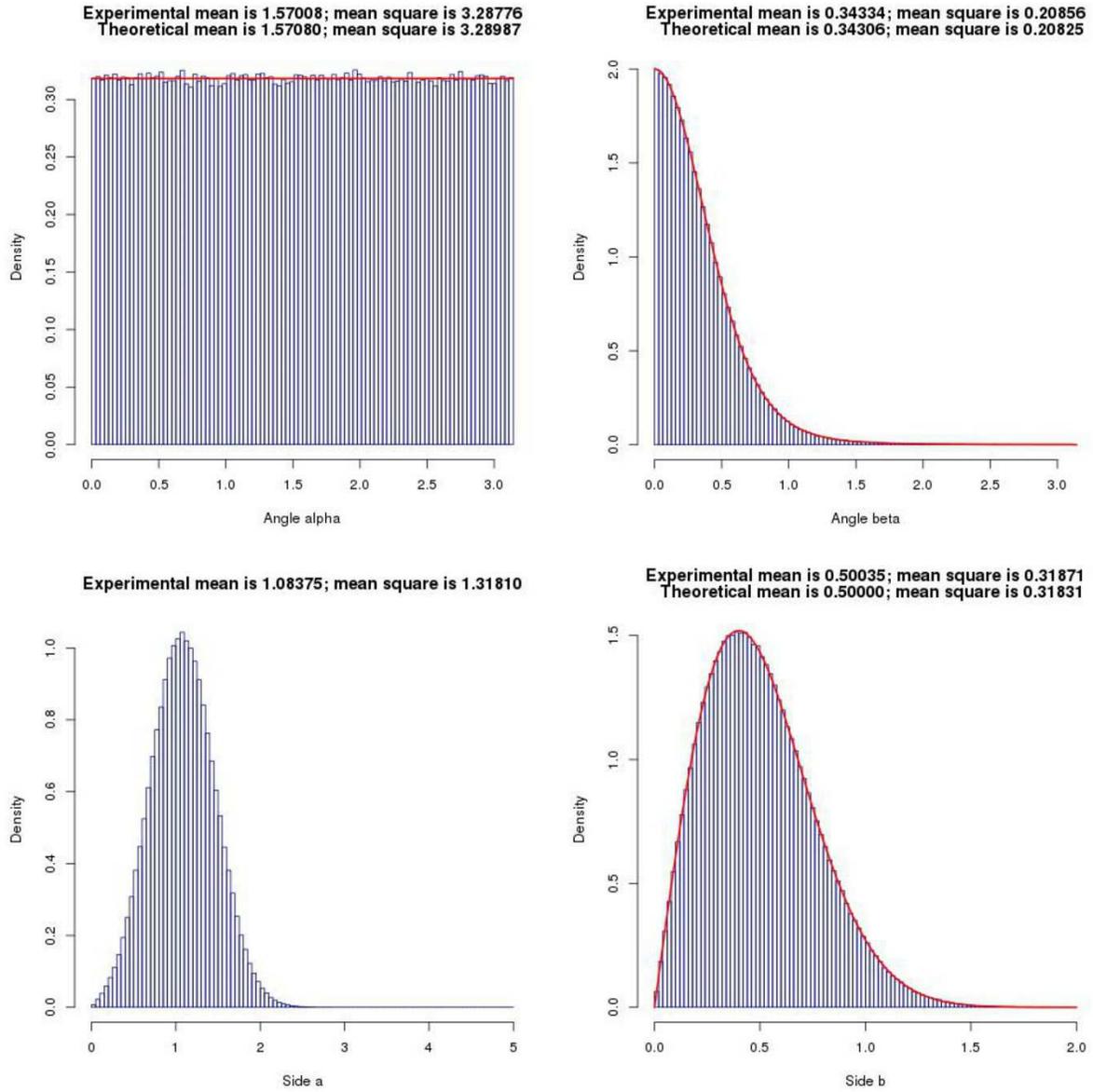}%
\caption{Staked densities for Poissonian triangle sides and angles.}%
\end{figure}

%

\begin{figure}[ptb]%
\centering
\includegraphics[
height=6.3362in,
width=6.3362in
]%
{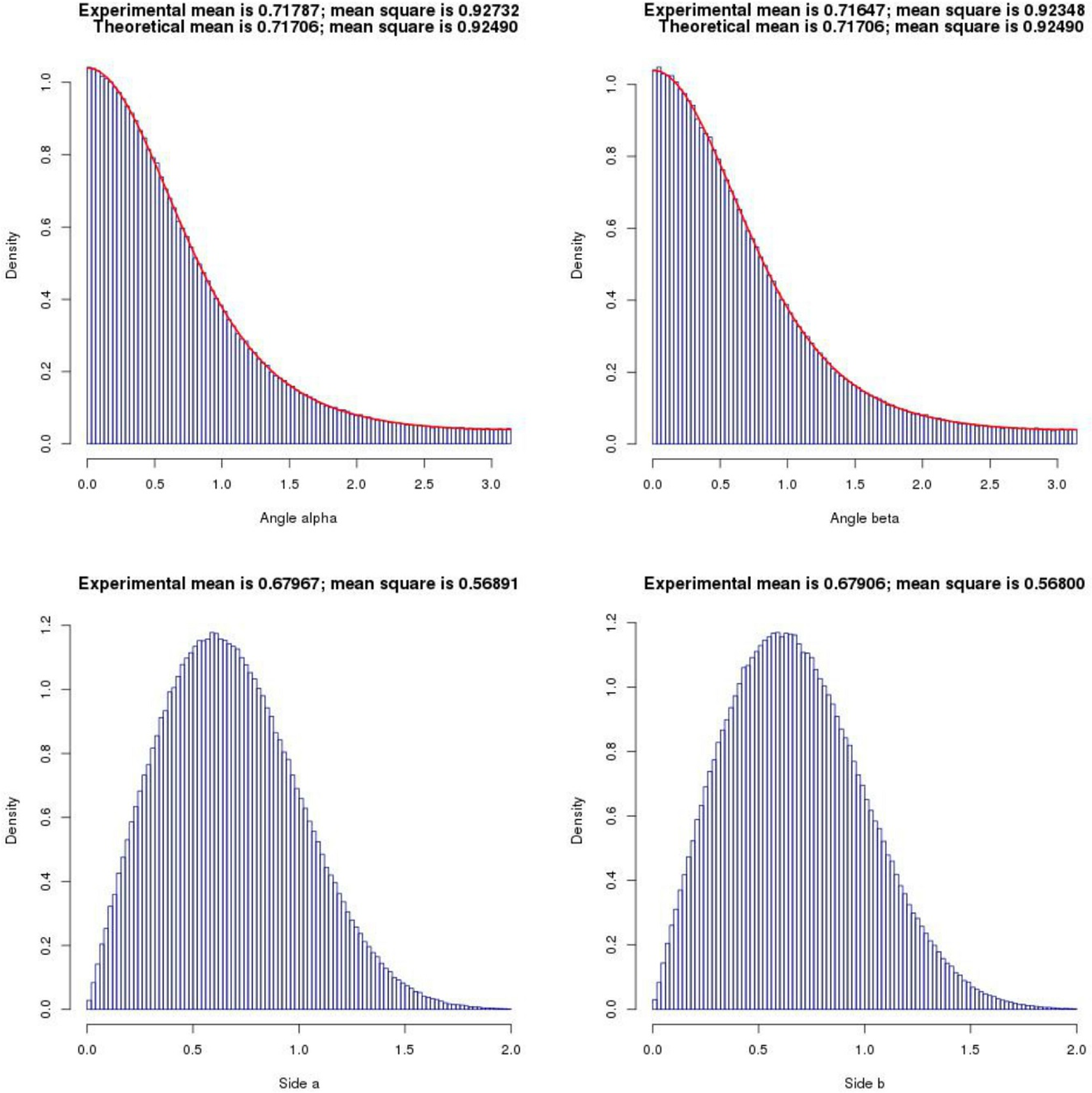}%
\caption{Anchored densities for Poissonian triangle sides and angles.}%
\end{figure}


\begin{thebibliography}{99}                                                                                               %


\bibitem {SS-Pois}D. Stoyan and H. Stoyan, \textit{Fractals, Random Shapes and
Point Fields: Methods of Geometrical Statistics}, Wiley, 1994, pp. 212--215;
MR1297125 (95h:60016).

\bibitem {Ha-Pois}F. Haken, \textit{Quantum Signatures of Chaos}, 2nd ed.,
Springer-Verlag, 2004, pp. 345--346, 388--389; MR2242927 (2008h:81055).

\bibitem {SN-Pois}J. Sakhr and J. M. Nieminen, Wigner surmises and the
two-dimensional homogeneous Poisson point process, \textit{Phys. Rev. E} 73
(2006) 047202.

\bibitem {ES1-Pois}B. Eisenberg and R. Sullivan, Random triangles in $n$
dimensions, \textit{Amer. Math. Monthly} 103 (1996) 308--318; MR1383668 (96m:60025).

\bibitem {Fi2-Pois}S. R. Finch, Random Gaussian tetrahedra, http://arxiv.org/abs/1005.1033.

\bibitem {DS-Pois}D. J. Daley, S. Ebert and R. J. Swift, Size distributions in
random triangles, \textit{J. Appl. Probab.} 51A (2014) 283--295; MR3317364.

\bibitem {Fi3-Pois}S. R. Finch, Catalan's constant, \textit{Mathematical
Constants}, Cambridge Univ. Press, 2003, pp. 53--59; MR2003519 (2004i:00001).

\bibitem {Fi4-Pois}S. R. Finch, Correlation between angle and side, http://arxiv.org/abs/1012.0781.

\bibitem {Fi5-Pois}S. R. Finch, Rank-$3$ projections of a $4$-cube, http://arxiv.org/abs/1204.3468.

\bibitem {Fi6-Pois}S. R. Finch, Pins, stakes, anchors and Gaussian triangles, http://arxiv.org/abs/1410.6742.

\bibitem {DL-Pois}D. J. Daley, S. Ebert and G. Last, Two lilypond systems of
finite line-segments, \textit{Probab. Math. Statist.} 36 (2016) 221--246;
http://arxiv.org/abs/1406.0096; MR3593022.

\bibitem {Grf-Pois}D. Griffiths, Uniform distributions and random triangles,
\textit{Math. Gazette}. 67 (1983) 38--42.

\bibitem {Mor-Pois}T. Moore, Random triangle problem (long summary), http://mathforum.org/kb/plaintext.jspa?messageID=86196.

\bibitem {ES2-Pois}A. Edelman and G. Strang, Random triangle theory with
geometry and applications, \textit{Found. Comput. Math.} 15 (2015) 681--713;
http://arxiv.org/abs/1501.03053; MR3348170.

\bibitem {BRT-Pois}A. Baddeley, E. Rubak and R. Turner, \textit{Spatial Point
Patterns: Methodology and Applications with R}, Chapman and Hall/CRC Press,
2015; http://spatstat.org/.

\bibitem {Fi7-Pois}S. R. Finch, Simulations in R\ involving triangles and
tetrahedra, unpublished software code (2017).

\bibitem {Fi1-Pois}S. R. Finch, Random triangles, unpublished note (2010).

\bibitem {Pa-Pois}A. Papoulis, \textit{Probability, Random Variables, and
Stochastic Processes}, McGraw-Hill, 1965, pp. 187--206; MR0176501 (31 \#773).
\end{thebibliography}
\end{document}